# EXACT SOLUTION OF DISCRETE HEDGING EQUATION FOR EUROPEAN OPTION

Yakovlev D.E[1]., Zhabin D. N.[2]

*Department of Higher Mathematics and Mathematical Physics*
*Tomsk Polytechnic University, Tomsk, Lenina avenue 30, 634004, Russia*

The approach that allows find European option price on the assumption of hedging at discrete times is proposed. The routine allows find the option price not for lognormal distribution functions of underlying asset only but for wide enough classes of distribution functions too. It is shown that there exists a nonzero possibility that market parameters can take values such that to realize the hedging policy becomes impossible. This fact is not in contradiction with Black-Scholes option price model as long as this possibility tends to zero at the limit of continuous hedging.

---
[1] *Email addresses: ENMFnauka@yandex.ru*
[2] *Email addresses: zhabin@phys.tsu.ru*

**Introduction**

At 1973 American scientists Black and Scholes (1973) proposed the pricing model for European option. At the time being Black-Scholes model is the commonly used one for pricing derivatives (see Hull (1999), Wilmott et al. (1996)). Note that this model is developing up to now (Korn and Korn (2000), Wilmott (1998, part 3)). The model is based on assumptions that all basically say three things: investor deals on efficient market; the underlying asset price follows the lognormal random walk; trading of underlying asset can take place continuously.

It should be noted that there are many faults with Black-Scholes model (e.g., Wilmott, 1998, section 19). Some of these are solved (e.g., Sircar and Papanicolaou (1998), Whalley and Wilmott (1997), Krakovsky (1999)). But instead of solved problems new ones are coming. For instance, as it follows from previously published works (e.g., Peters (1994)) the hypothesis of assets lognormal random walk is doubted. Moreover, in Wilmott (1998, section 20) is shown that a generalization of Black-Scholes model to the hedging at discrete times is actual since continuous hedging is impossible, and even undesirable, in practice.

In this paper we propose the extension of Black-Scholes model to wide enough classes of distribution functions of assets return in assumption that trading of the underlying asset can take place at discrete times only. In contrast to Whalley and Wilmott (1997) we present the exact way how to find the option price without any assumption about negligibility of the contribution caused by discrete hedging. We found that as the hedging periods tend to zero and the numbers of hedging acts comes to infinity the result of Whalley and Wilmott (1997) can be found with prescribed accuracy.

**Approach**

In the approach proposed we assume that the underlying asset price is the Markov process (see Gramer and Leadbetter, 1967). The asset price $\{S_t^{market}\}_{t \geq t_0}$ changes together with a time and the time is continual. Each cross-section $S_t^{market}$ is a stochastic variable that represents the asset price at time $t$. Since in a real market asset prices are quoted at discrete time intervals we assume that an investor deals with period $\tau$.

By definition, put $t_k = T - k\tau$, $k = 1, 2, \ldots, n$, where $n$ is defined from equality $n\tau = T - t_0$. Here $T$ is an option maturity; $t_0$ is a "starting date". This means that for the investor the asset price is represented by a sequence of stochastic variables $\{S_k\}_{k \geq k_0}$, $S_k = S_{t_k}^{market} = S_{T-k\tau}^{market}$. Without loss of generality we will assume further that $t_0 = 0$. In addition, we consider auxiliary variables $\xi_k$ such that

$$\xi_k : \xi_k = \ln\left(\frac{S_k}{S_{k+1}}\right), \text{ i.e. } S_k = e^{\xi_k} S_{k+1}. \tag{1}$$

To be define, we assume that for $\forall p \in (-a, a)$, where $a > 2$ there exist expectations of $e^{p \cdot \xi_k}$ such that

$$E[e^{p \cdot \xi_k}] < \infty. \tag{2}$$

Further, at time $t_{k+1}$ we construct a portfolio consisting of one long option position and a short position in some quantity $\Delta_k$ of the underlying. Let us use $\Pi_{k+1}$ to define the value of this portfolio at time $t_{k+1}$. By construction, we get

$$\Pi_{k+1}(s) = V_{k+1}(s) - \Delta_{k+1} \cdot S_{k+1}. \tag{3}$$

Now assume that one time-step was passed. It is evident that $\Delta_{k+1}$ is held fixed during the time-step. Let us denote by $S_{k+1} = s$ the asset price at time $t_{k+1}$ and by $D[A]$ a variance of the stochastic variable $A$. Then for the portfolio value at time $t_k$ we get

$$\widetilde{\Pi}_k = V_k(S_k) - \Delta_{k+1} \cdot S_k = V_k(e^{\xi_k} \cdot s) - \Delta_{k+1} \cdot s \cdot e^{\xi_k}. \tag{4}$$

By carefully choosing $\Delta_{k+1}$ we can eliminate the variance of portfolio (4). It can be shown in the usual way that as soon as $\Delta_{k+1}$ is equal to

$$\Delta_{k+1} = \frac{1}{s} \frac{\text{cov}[V_k(se^{\xi_k}), e^{\xi_k}]}{D[e^{\xi_k}]}$$ the portfolio variance becomes minimal.

Let us note that Black-Scholes pricing policy may be written as

$$E[\widetilde{\Pi}_k] = e^{r_k \cdot \tau} \cdot \Pi_{k+1}. \tag{5}$$

Here $r_k$ is a risk–free interest rate at time $t_k$. It is easy to check that equation (5) is equal to the integral equation

$$\int_{R^1} V_k(se^x) f_k(x) dx = V_{k+1}(s). \tag{6}$$

Here $f_k(x) = \left( \dfrac{(e^x - m_k)(1 - e^{-r_k \tau} m_k)}{d_k} + e^{-r_k \tau} \right) \cdot u_k(x)$; $u_k(x)$ is the distribution density of quantities $\xi_k$; $m_k = E[e^{\xi_k}]$; $d_k = D[e^{\xi_k}]$. Since the function $V_0(s)$ is assumed to be known we may express $V_k(s)$ for whole time sequence.

Let us remark that at expiry the value of the call option (payoff) can be written as $V_0(s) = \max[s - E, 0]$. Taking into account the inequality $|V_0(s)| \leq \dfrac{s^m}{m \cdot E^{m-1}}$ for $1 \leq m < a - 1$ we can prove that there exist numbers $A_k(m)$ such that $|V_k(s)| \leq A_k(m) \cdot s^m$. In the same way it can be shown that all $V_k(s)$ are continuous and there exist numbers $W_k$ such that the inequality $|V_k(s) - s| \leq W_k$ is true. The last inequality means that $\dfrac{V_k(s)}{s} \approx 1$ for $s \gg E$.

To solve (6), let us rewrite this equation with help of Mellin transformation. Recall that Mellin transform of a function $h(x)$ is the integral

$$H(p) = \int_0^\infty x^{p-1} h(x) dx, \qquad (7)$$

if this integral exist (see Doetsch (1954), Churchill (1956)). Further, by $H(p) = Mel[h(x)]$ we will denote Mellin transform; by $h(x) = Mel^{-1}[H(p)]$ we will denote the inverse Mellin transform. It is clear that the domain of the function $H(p)$ is the set of complex $p$ such that integral (7) converges.

From properties of functions $V_k(s)$ discussed above it follows that for any $p$ such that $1 - a < \operatorname{Re} p < -1$ there exist integrals

$$F_k(p) = \int_0^\infty x^{p-1} V_k(x) dx.$$

Moreover, from (2) it follows that the function $f_k(\ln y)$ admits Mellin transformation for $1 - a \leq \operatorname{Re} p \leq a - 1$. By definition, put

$$U_k(p) = Mel[f_k(\ln x)]. \tag{8}$$

It follows that equation (6) can be written with help of Mellin transformation in form

$$F_{k+1}(p) = U_k(-p) \cdot F_k(p). \tag{9}$$

Obviously, the solution of this difference equation is given by the following expression

$$F_k(p) = F_0(p) \cdot \prod_{m=0}^{k-1} U_m(-p). \tag{10}$$

Here $F_0(p)$ is Mellin transform of the payoff function. Thus we found the solution of equation (6) in term of Mellin transform. There are two alternative ways to reconstruct functions $V_k(s)$. In the first of place we can use the inverse Mellin transformation. Therefore, we get

$$h(x) = Mel^{-1}[H(p)] = \frac{1}{2\pi i} \int_{a-i\infty}^{a+i\infty} H(p) x^{-p} dp. \tag{11}$$

In this way we obviously obtain

$$V_k(s) = \frac{1}{2\pi i} \int_{a_0 - i\cdot\infty}^{a_0 + i\cdot\infty} F_k(p) \cdot s^{-p} \cdot dp, \text{ where } a_0 = \frac{a}{2}. \tag{12}$$

In the second place we define auxiliary functions

$$G_k(x) = Mel^{-1}\left[\prod_{m=0}^{k-1} U_m(-p)\right]. \tag{13}$$

Using properties of Mellin transformation, it follows that

$$V_k(s) = \int_0^\infty G_k\left(\frac{s}{y}\right) \cdot \frac{1}{y} \cdot V_0(y) dy. \tag{14}$$

Here function $G(k,x,y) = G_k\left(\frac{x}{y}\right) \cdot \frac{1}{y}$ is Green function.

To be precise, $\tau$ is usually small by value with respect to the time till maturity. Thus it is useful to build an asymptotic expansion with respect to the small $\tau$. As soon as $n\tau = T$ we can denote by $t = k\tau$ any time till maturity, here $k=1, 2, \ldots, n$. This means that $V_k(s) = V_{\left[\frac{t}{\tau}\right]}(s) = V(t,s)$. Further, expand $V(t,s)$ into a series in powers of $\tau$. Therefore, we get

$$V(t,s) \approx B_0(t,s) + \tau B_1(t,s) + \tau^2 B_2(t,s) + \ldots + \tau^l B_l(t,s). \tag{15}$$

Expression (15) allows estimate the option price at time $t$ for small hedging period. It can easily be checked that in case of the lognormal random walk of the underlying the function $V(t,s)$ approximates the solution of Whalley and Wilmott (1997) with any prescribed accuracy.

Note also that from analysis of expression (6) it is follows that for $k \geq 1$ there exist nonempty sets $\Xi_k \subset [0,+\infty)$ such that the function $V_k(s)$ becomes negative for $s \in \Xi_k$. In other words this means that there exist possible states of the market such that to realize the hedging policy becomes impossible. To avoid this problem we come to restriction on model parameters

$$\frac{1}{\tau} \ln\left(E[e^{\xi_k^\tau}]\right) \leq r_k \leq \frac{1}{\tau} \ln\left(\frac{E[e^{2\xi_k^\tau}]}{E[e^{\xi_k^\tau}]}\right). \tag{16}$$

Since $D[e^{\xi_k^\tau}] = E[e^{2\xi_k^\tau}] - E^2[e^{\xi_k^\tau}] > 0$, it follows that the interval (16) is not empty. Moreover, as long as we come to the continuous hedging, i.e. $\tau \to 0+$ while $k\tau = T$, sets $\Xi_k$ degenerates to empty sets. If we suppose that the underlying asset is the lognormal distributed stochastic variable, then the relationship between Black-Scholes option price $C^{BS}(s,t)$ and option price $V(t,s)$ is true

$$C^{BS}(s,t) = \lim_{\tau \to +0} V(t,s). \qquad (17)$$

**Example**

Let us consider the example of the exact expression for the option price with hedging at discrete times. Suppose that for $\forall k$, $\xi_k = \xi$, where $\xi$ is the normally distributed stochastic variable with mean $E[\xi] = \mu\tau$ and variance $D[\xi] = \sigma^2\tau$. In this case, the probability density function is given by

$$u_k(x) = u(x) = \frac{1}{\sigma\sqrt{2\pi\tau}} e^{-\frac{(x-\mu\tau)^2}{2\sigma^2}}.$$

Evidently, we have

$$m_k = m = e^{\mu\tau + \frac{\sigma^2\tau}{2}}, \quad d_k = d = e^{2\mu\tau + \sigma^2\tau}\left(e^{\sigma^2\tau} - 1\right).$$

We assume also that all $r_k$ are constant till option maturity, i.e. $r_k = r$. Combining all together, we obtain

$$U_k(p) = U(p) = e^{\mu\tau p + \frac{\sigma^2\tau}{2}p^2}\left(M_1 + M_2 e^{p\sigma^2\tau}\right).$$

Here $M_1$ and $M_2$ are given by

$$M_1 = e^{-r\tau} + \frac{e^{(\mu-r)\tau+\frac{\sigma^2}{2}\tau} - 1}{e^{\mu\tau+\frac{3}{2}\sigma^2\tau} - e^{\mu\tau+\frac{\sigma^2\tau}{2}}} \quad , \quad M_2 = \frac{1 - e^{(\mu-r)\tau+\frac{\sigma^2}{2}\tau}}{e^{\mu\tau+\frac{3}{2}\sigma^2\tau} - e^{\mu\tau+\frac{\sigma^2\tau}{2}}}.$$

So, the solution of equation (10) can be written in form

$$F_n(p) = F_0(p)U^n(-p) = F_0(p) \cdot e^{-n\mu\tau p + n\frac{\sigma^2\tau}{2}p^2} \left( M_2 e^{-\sigma^2\tau p} + M_1 \right)^n =$$

$$= F_0(p) \sum_{k=0}^{n} C_n^k M_1^{n-k} M_2^k e^{-p(n\mu\tau+k\sigma^2\tau)+n\frac{\sigma^2\tau}{2}p^2} = F_0(p) \cdot L_n(p).$$

Hence, we obtain

$$L_n(p) \xrightarrow{Mel} \sum_{k=0}^{n} C_n^k M_1^{n-k} M_2^k \frac{1}{\sigma\sqrt{2\pi n\tau}} e^{-\frac{(\ln x + n\mu\tau + k\sigma^2\tau)^2}{2n\sigma^2\tau}} = G_n(x).$$

It follows that

$$V_n(s) = \int_0^\infty G_n\left(\frac{s}{x}\right)\frac{1}{x} V_0(x) dx.$$

Changing the variable $y = \frac{1}{x}$, we get

$$V_n(s) = \sum_{k=0}^{n} C_n^k M_1^{n-k} M_2^k \frac{1}{\sigma\sqrt{2\pi n\tau}} \int_0^{\frac{s}{E}} e^{-\frac{(\ln y + n\mu\tau + k\sigma^2\tau)^2}{2n\sigma^2\tau}} \left(\frac{s}{y} - E\right) dy. \quad (18)$$

It can be shown in standard way that for continuous hedging this expression of the option price is equal to the well known Black-Scholes formula.

**Summary**

In present paper we propose the method that allows find the exact expression for price of European option on assumption of hedging at discrete times. Another advance of the method is that this method provides the tool to price option for wider classes of underlying asset density function then lognormal density functions.

We show that the approach gives the well known Black-Scholes formula of option price at the limit of continuous hedging and on the assumption of lognormal density function of the underlying assets. Moreover, we consider the asymptotic behavior of option price for small hedging period. It is shown that this asymptotic coincides with result Whalley and Wilmott (1997) with accurate to the first order.

Finally, we show that there exist some market states such that to realize hedging becomes impossible. This fact is not in contradiction with Black-Scholes theory as long as we come to the continuous hedging sets $\Xi_k$ degenerates to empty sets.